\documentclass[12pt]{amsart}

\usepackage{a4,amsmath,amssymb,amsfonts}
\usepackage[all]{xy}

\newtheorem{theorem}{Theorem}[section]
\newtheorem{lemma}[theorem]{Lemma}
\newtheorem{proposition}[theorem]{Proposition}
\newtheorem{corollary}[theorem]{Corollary}

\begin{document}
\baselineskip=15.5pt

\renewcommand{\phi}{\varphi}

\newcommand{\Gal}{{\rm Gal}}
\newcommand{\Grass}{{\rm Grass}}
\newcommand{\Mat}{{\rm Mat}}
\newcommand{\Ocal}{{\mathcal O}}
\newcommand{\PGL}{{\rm PGL}}
\newcommand{\rarpa}[1]{\stackrel{#1}{\longrightarrow}}
\newcommand{\rk}{{\rm rk}}
\newcommand{\Spec}{{\rm Spec}}

\newcommand{\adop}{{\mathbb A}}
\newcommand{\cdop}{{\mathbb C}}
\newcommand{\ndop}{{\mathbb N}}
\newcommand{\pdop}{{\mathbb P}}
\newcommand{\zdop}{{\mathbb Z}}

\title{Generalization of a criterion for semistable vector bundles}

\author[I. Biswas]{Indranil Biswas}

\address{School of Mathematics, Tata Institute of Fundamental
Research, Homi Bhabha Road, Bombay 400005, India}

\email{indranil@math.tifr.res.in}

\author[G. Hein]{Georg Hein}

\address{Universit\"at Duisburg-Essen, Fachbereich
Mathematik, 45117 Essen, Germany}

\email{georg.hein@uni-due.de}

\begin{abstract}
It is known that a vector bundle $E$ on a smooth projective curve
$Y$ defined over an algebraically closed field is semistable if
and only if there is a vector bundle $F$ on $Y$ such that both
$H^0(X,\, E \otimes F)$ and $H^1(X,\, E \otimes F)$ vanishes.
We extend this criterion for semistability to vector bundles 
on curves defined over perfect fields.
Let $X$ be a geometrically irreducible smooth projective curve
defined over a perfect field $k$, and let $E$ be a vector bundle
on $X$. We prove that $E$ is semistable if and only if
there is a vector bundle $F$ on $X$ such that
$H^i(X,\, E\otimes F)\, =\, 0$ for all $i$. We also give an
explicit bound for the rank of $F$.
\end{abstract}

\maketitle

\section{Introduction}

A theorem due to Faltings says that
a vector bundle $E$ on a smooth projective curve
$Y$ defined over an algebraically closed field of characteristic
zero is semistable if
and only if there is a vector bundle $F$ on $Y$ such that both
$H^i(X,\, E\otimes F)\, =\, 0$ for all $i$ \cite[p. 514, Theorem 
1.2]{Fa}. It is is known that this criterion for semistability
extends to vector bundles on smooth projective curves defined
over an algebraically closed fields of positive characteristic.
(See \cite{He}, \cite{BH} for related results.)

Our aim here is to investigate this criterion for curves
defined over finite fields and more generally over perfect fields.
We prove the following theorem.

\begin{theorem}\label{thm1}
Let $X$ be a geometrically irreducible smooth projective curve
defined over a perfect field $k$.
A vector bundle $E$ over $X$ is semistable if and only if there
is a vector bundle $F$ over $X$ such that
$H^i(X,\, E\otimes F)\, =\, 0$ for all $i$.
\end{theorem}

We also produce an effective bound for the rank of $F$ in
Theorem \ref{thm1}. More precisely, given nonnegative
integers $g$ and $r$, and an integer $d$, there is an
explicit integer
$R(g,r,d)$ such that for any triple $(k\, ,X\, ,E)$, where
\begin{itemize}
\item $k$ is a perfect field,

\item $X$ is a geometrically irreducible smooth projective curve
of genus $g$ defined over $k$, and

\item $E$ is a vector bundle over $X$ of rank $r$ and degree $d$,
\end{itemize}
the vector bundle $E$ is semistable if and only
if there is a vector bundle $F$
over $X$ of rank $R(g,r,d)$ such that $H^i(X,\,
E\otimes F)\, =\, 0$ for all $i$. (See Theorem \ref{crit}.)

Let $X$ be a scheme defined over a field $k$, and 
let $D\,\subset\, X$ be an effective divisor. Then there is a
geometric point of $X$ that lies outside $D$.
Corollary \ref{ed} bounds the degree of
the field extension $K/k$ such that $D(K) \,\subsetneqq\, X(K)$.
This is used in the proof of Theorem \ref{thm1}.

\section{Rational points outside a given
hypersurface}\label{EffDiv}

Let $k$ be any field. The algebraic closure of $k$ will
be denoted by $\overline k$.

\begin{lemma}\label{effdiv0}
Let $D \,\subset\, \adop^n_k$ be an effective divisor defined
over $\overline k$. Given any field extension $K/k$
such $K$ has more than $\deg(D)$ elements, there exists a
$K$-rational point in $\adop^n_K$ that lies outside $D$.
\end{lemma}

\begin{proof}
One follows the proof of Proposition 1.3(a) in
\cite[p. 4]{Kunz} almost
word for word simply replacing ``{\em infinite field}'' by 
``{\em field with more than $\deg(D)$
elements}'': We assume that $D$ is given by the polynomial $F\,\in
\, \overline{k}[X_1,\ldots,X_n]$, and proceed by induction over $n$.
For $n\,=\,1$ it is the statement that a polynomial of degree $d$
cannot have more than $d$ zeros.
Now assume that $X_n$ occurs in $F$ and write
$$
F\,= \, \phi_0 + \phi_1X_n + \ldots \phi_tX_n^t\, ,
$$
where $\phi_i \,\in \, \overline{k} [X_1,\ldots,X_{n-1}]$
and $\phi_t \, \ne\, 0$.
Since $t$ and $\deg(\phi_t)$ are both at most $\deg(D)$, we
conclude from the induction hypothesis
the existence of a point $(x_1,\ldots ,x_{n-1}) \,\in\,
K^{n-1}$ such that $\phi_t(x_1,\ldots ,x_{n-1}) \,\ne\, 0$.
Now the polynomial $X_n \,\longmapsto\,
F(x_1,\ldots ,x_{n-1},X_n)$ has at most $t$ zeros.
\end{proof}

\begin{lemma}\label{effdiv1}
Let $D \,\subset\, \pdop_k^n$ be an effective divisor in
projective space $\pdop^n$ defined over $\overline{k}$. Then for
any extension $K/k$ such that $K$ has at least $\deg(D)$
elements, there is a $K$-rational point in $\pdop^n(K)$
that lies outside $D(K)$.
\end{lemma}

\begin{proof}
Again we proceed by induction on $n$. The case $n\,=\,1$ is obvious.
For $n\,>\,1$, we consider the pencil of hyperplanes defined over $K$
passing through a
codimension two linear subspace. Since there are more than $\deg(D)$
of these hyperplanes, the union of all these hyperplanes
cannot be contained in $D$. Thus, there exists a hyperplane
$H \,\cong \,\pdop^{n-1} \,\subset\, \pdop^n$ that intersects $D$ 
properly. Now the proof is completed by the induction hypothesis.
\end{proof}

Let $\Grass(m,n)$ be the Grassmannian of $m$--dimensional linear
subspaces of $k^n$. Let
\begin{equation}\label{pluk}
\iota\, :\,\Grass(m,n)\, \longrightarrow\,
{\mathbb P}\, :=\, \pdop^{\binom{n}{m} -1}
\end{equation}
be the Pl\"ucker embedding. By an hypersurface of degree $d$
on $\Grass(m,n)$ we will mean one from the complete linear
system $\vert \iota^*{\mathcal O}_{\mathbb P}(d)\vert$.

\begin{lemma}\label{effdiv2}
Let $D \,\subset \,\Grass(m,n)$ be a hypersurface
in the Grassmannian. If a field extension $K/k$ has more than
$m\cdot \deg(D)$ elements, then there is a $K$--rational point of
$\Grass(m,n)(K)$ that is not contained in $D(K)$.
\end{lemma}

\begin{proof}
We consider the dense open cell in the Grassmannian given by the
open immersion $j:\adop^{m\cdot(n-m)}\, \longrightarrow\,
\Grass(m,n)$ defined by
\[ \left( a_{i,j} \right)_{i=1,\ldots,m \,\,\, j:=m+1,\ldots,n}
\longmapsto
{\rm span}\left( \begin{array}{cccccccc}
1&0& \cdots & 0 & a_{1,m+1} & a_{1,m+2} & \cdots & a_{1,n} \\
0&1& \cdots & 0 & a_{2,m+1} & a_{2,m+2} & \cdots & a_{2,n} \\
\vdots&\vdots & \ddots & \vdots & \vdots &\vdots& \ddots & \vdots\\
0&0& \cdots & 1 & a_{m,m+1} & a_{m,m+2} & \cdots & a_{m,n} \\
\end{array}
\right) \]
The Pl\"ucker embedding (see \eqref{pluk})
restricted to $\adop^{m\cdot(n-m)}$ is given by
the $m \times m$--minors of degree at most
$m$ of the above matrix. Therefore,
$j^*D$ is a divisor of degree $m \cdot \deg(D)$.
Now the proof is completed using Lemma \ref{effdiv0}.
\end{proof}

\begin{proposition}\label{effdiv}
Let ${\mathcal O}_X(H)$ be a globally generated ample line bundle
on a projective scheme $X$ of dimension $n$ defined over $k$. 
Let $D \,\subset\, X$ be an effective divisor $D \subset X$. Let
$K_1/k$ be a field extension that has more than $\max\{ (n+1)H^n, 
D.H^{n-1}-1\}$ elements. Then there exists a
field extension $K_2/K_1$ with $[K_2:K_1] \leq H^n$,
such that there is a $K_2$--rational point of $X(K_2)$
that does not lie in $D(K_2)$.
\end{proposition}
\begin{proof}
We consider the short exact sequence of vector bundles
$$
0\, \longrightarrow\, W \, \longrightarrow\, H^0(L) \otimes
\Ocal_X\, \longrightarrow\, L \, \longrightarrow\, 0
$$
over $X$. Let $\Grass_X(n+1,W)$ be the Grassmann bundle over $X$
parameterizing all $(n+1)$--dimensional
subspaces in the fibers of $W$. We have
\[ \begin{array}{rcl}
\dim \Grass_X(n+1,W)&=&\dim(X)+(h^0(L)-n-2)(n+1)\\
&=&  \dim(\Grass(n+1,H^0(L)))-1\, .\\
\end{array}\]

We will show that the degree of the hypersurface
$\Grass_X(n+1,W)$ in $\Grass(n+1,H^0(L))$ is $H^n$. 
To prove this, take any subspace $U \subset H^0(X, L)$ of
dimension $n+2$.
The $(n+1)$-dimensional subspaces of $U$ form a
projective line $\pdop^1_k$ in $\Grass(n+1,H^0(L))$. The
degree of the restriction of $\iota^*{\mathcal O}_{\mathbb P}(1)$
(see \eqref{pluk}) to this $\pdop^1_k$
is one. To compute the intersection number of the line with
$\Grass_X(n+1,W)$ we may assume that $H^0(X,\, L)\,=\,U$.
So it suffices to count the intersection of
a line in $\pdop(U)$ with the
divisor $X \,\subset\, \pdop(U)$. Thus, we conclude that
the hypersurface $\Grass_X(n+1,W)\, \subset\,
\Grass(n+1,H^0(L))$ is of degree $H^n$.

Now, using Lemma \ref{effdiv2} and the assumption on $K_1$ we
conclude that there exits a $K$--point in $\Grass(n+1,H^0(X, L))$
not lying in $\Grass_X(n+1,W)$. This yields a finite morphism $X 
\rarpa{\pi} \pdop^n$ defined over $K_1$. Now $\pi_*(D)$ is a
divisor of degree $D.H^{-1}$ on $\pdop^n$. Our assumption on the
number of elements in $K_1$ and Lemma \ref{effdiv1} together
imply that there is a $K_1$--rational point $P$ in the
complement of $\pi_*(D)$ in $\pdop^n$. The morphism $X_P 
\,\longrightarrow\, \Spec(K_1)$ is finite of degree $H^n$, and
it is defined over $K_1$. Thus, we 
find at least one point in $X_P$ defined over a field $K_2$ as in 
the statement of the proposition. This completes the proof of
the proposition.
\end{proof}

Proposition \ref{effdiv} has the following corollary.

\begin{corollary}\label{ed}
Given positive integers $n$, $\alpha$ and $\beta$,
define
$$
M(n,\alpha,\beta)\, :=\, \alpha \lceil \log_2(
\max\{ (n+1)\alpha+1,\beta \}) \rceil\, .
$$
For any quadruple $(k\, ,X\, ,H\, ,D)$, where
\begin{itemize}
\item $k$ is a field,

\item $X$ is a projective scheme of dimension $n$ defined over $k$,

\item $H\, \subset\, X$ is a base point--free ample hypersurface
with $H^n\, =\,\alpha$, and

\item $D\, \subset\, H$ is an effective divisor with $D.H^{n-1}\,
=\, \beta$,
\end{itemize}
there is a field extension $K/k$ of degree
$[K:k] \, \leq\, M(n,\alpha,\beta)$ with the property that
$X(K)$ has a $K$--rational point that does not lie in $D(K)$.
\end{corollary}

If we restrict ourselves only to infinite fields, then
$M(n,\alpha,\beta)$ in Corollary \ref{ed} can be taken to
be $\alpha$. If we fix a prime $p$ and
restrict ourselves only to fields of characteristic $p$, then
$M(n,\alpha,\beta)$ in Corollary \ref{ed} can be taken to
be $\alpha \lceil \log_p( \max\{ (n+1)\alpha +1,\beta \}) 
\rceil$.

\section{Semistability criterion over perfect fields}
\label{Crit}

\begin{theorem}\label{crit}
Let $X$ be a geometrically irreducible smooth projective curve
of genus $g$ defined over a perfect field $k$. Fix a positive
integer $r$ and an integer $d$. Then there is an explicit positive 
integer $R$ that depends only on $r$, $d$ and $g$ (in particular,
$R$ is independent of $k$) with the following property: A
vector bundle $E$ over $X$ of rank $r$ and degree $d$ is
semistable if and only if there is a vector bundle $F$ over $X$
of rank $R$ such that $H^i(X,\,E \otimes F)\,=\,0$ for all $i$.
\end{theorem}

\begin{proof}
If $E$ is not semistable, then clearly there is no $F$
such that $H^i(X,\,E \otimes F)\,=\,0$ for all $i$.
Let $E$ be a semistable vector bundle over $X$ of rank $r$
and degree $d$. We will construct $R$ and $F$.

The moduli space of semistable vector bundles over $X$ of rank
$r'$ and degree $d'$ will be denoted by ${\mathcal U}_X(r',d')$.

Let $h\,:=\,\gcd(r,d)$. Furthermore, we set
$\overline{r}\,:=\,\frac{r}{h}$, 
and $\overline{d}\,:=\, \frac{d}{h}$.
For any integer $n \geq 1$, consider the morphism
${\mathcal U}_X(n\overline{r},n (\overline{r}(g-1)-
\overline{d}))\, \longrightarrow\, 
{\mathcal U}_X(nr\overline{r}, nr\overline{r}(g-1))$
defined by $V\, \longmapsto\, V\otimes E$. Let
$\Theta_E$ denote the pull back of the natural theta
divisor in ${\mathcal U}_X(nr\overline{r}, nr\overline{r}(g-1))$
by this morphism. Therefore, $\Theta_E\times_k\overline{k}$ 
consists of all semistable vector bundles $W$ over
$X_{\overline{k}}\, =\, X\times_k \overline{k}$
of rank $n\overline{r}$ and degree $n(\overline{r}(g-1)-
\overline{d})$ such that $H^0(X_{\overline{k}},\,
W\otimes (E\otimes_k \overline{k}))\, \not=\, 0$.
(We note that a vector bundle $V'$ over $X$ is semistable
if and only if the vector bundle $V'\otimes_k \overline{k}$
over $X_{\overline{k}}$ is semistable;
see \cite[p. 222]{HN}.) The subscheme
$\Theta_E$ defined above is either an effective Cartier divisor
in the complete linear system $| h \cdot \Theta |$ or it is
the entire moduli space ${\mathcal U}_X(n\overline{r},n
(\overline{r} (g-1)- \overline{d}))$ (cf. \cite[\S~0.2.1]{DN}).

Popa showed that this is indeed a divisor in the linear
system $| h \cdot \Theta |$ for all $n \,\geq\, \frac{r^2+1}{4}$
(see \cite[p. 490, Theorem 5.3]{Pop}). Let $n$ be the
smallest integer such that
$n\,\geq\, \frac{r^2+1}{4}$. Consider the
effective divisor $\Theta_E\,\subset\, {\mathcal U}_X( n
\overline{r},n (\overline{r} (g-1)- \overline{d}))$ for $n:=\lceil 
\frac{r^2+1}{4} 
\rceil$.
By Corollary \ref{ed}, there exists an integer $M$
and a field extension $K/k$ such that 
$\Theta_E$ does not contain all $K$--rational points of 
${\mathcal U}_X(n\overline{r},n
(\overline{r} (g-1)- \overline{d}))$.

The integer $R$ in the statement of the theorem will be
$n \overline{r} M!$.

Since $\Theta_E$ does not contain all $K$--rational points of
${\mathcal U}_X(n\overline{r},n
(\overline{r} (g-1)- \overline{d}))$, there exists a
vector bundle $F_1$ of rank $n \overline{r}$ defined over
$X_K\,=\, X\times_k K$, where $K/k$ is some Galois 
extension of degree dividing $M!$, such that
\begin{equation}\label{e1}
H^0(X_K,\, (E\otimes_k K) \otimes F_1)\,=\, 0
\,=\, H^1(X_K,\, (E\otimes_k K) \otimes F_1)\, .
\end{equation}
{}From \eqref{e1} it follows that
\begin{equation}\label{e2}
H^0(X_K,\, (E\otimes_k K) \otimes \sigma^*F_1)\,=\, 0
\,=\, H^1(X_K,\, (E\otimes_k K) \otimes \sigma^*F_1)
\end{equation}
for all $\sigma\, \in\, \Gal(K/k)$.

Now we consider the direct sum
$$
F_2\, := \, \bigoplus_{\sigma \in \Gal(K/k)}
\sigma^*F_1\, .
$$
This $F_2$ is a vector bundle defined over $X$. From \eqref{e2}
it follows immediately that
that $H^i(X, \, E \otimes F_2)=0$ for all $i$.
Also, the rank of $F_2$ clearly divides $n\overline{r} M!$.
Finally we set $m\,:=\, \frac{n \overline{r} M!}{\rk(F_2)}$,
and $F\,:= \,F_2^{\oplus m}$, and obtain the asserted vector
bundle of rank $R=n \overline{r} M!$.
\end{proof}

\end{document}